\documentclass[11pt,a4paper]{amsart}

\newcommand*\patchAmsMathEnvironmentForLineno[1]{%
  \expandafter\let\csname old#1\expandafter\endcsname\csname #1\endcsname
  \expandafter\let\csname oldend#1\expandafter\endcsname\csname end#1\endcsname
  \renewenvironment{#1}%
     {\linenomath\csname old#1\endcsname}%
     {\csname oldend#1\endcsname\endlinenomath}}%
\newcommand*\patchBothAmsMathEnvironmentsForLineno[1]{%
  \patchAmsMathEnvironmentForLineno{#1}%
  \patchAmsMathEnvironmentForLineno{#1*}}%
\AtBeginDocument{%
\patchBothAmsMathEnvironmentsForLineno{equation}%
\patchBothAmsMathEnvironmentsForLineno{align}%
\patchBothAmsMathEnvironmentsForLineno{flalign}%
\patchBothAmsMathEnvironmentsForLineno{alignat}%
\patchBothAmsMathEnvironmentsForLineno{gather}%
\patchBothAmsMathEnvironmentsForLineno{multline}%
\patchBothAmsMathEnvironmentsForLineno{eqnarray}%
}

\usepackage{amsthm, amsfonts, amssymb, amsmath, latexsym, enumerate, array}
\usepackage[all]{xy}
\usepackage{hyperref,mdwlist,mathrsfs}
\usepackage[left=2.8cm,top=2.7cm,right=2.7cm]{geometry}
\usepackage{lineno}

\newtheorem{thm}{Theorem}[section] \newtheorem{lem}[thm]{Lemma}
\newtheorem{corol}[thm]{Corollary} \newtheorem{prop}[thm]{Proposition}

\newtheorem*{thm*}{Theorem}

\theoremstyle{definition} 
 
\newtheorem{egs}[thm]{Examples}

\newtheorem{rem}[thm]{Remark}
\newenvironment{rem*}{\begin{rem}\em}{\end{rem}}
\newcommand {\GG}{\mathbb{G}}
\newcommand{\ra}{\rightarrow}
\newcommand{\TC}{\mathbb{T}}
\newcommand {\CC}{\mathbb{C}}
\newcommand{\cG}{\mathcal{G}}
\newtheorem{facts}[thm]{Facts}
\newenvironment{facts*}{\begin{facts}\em}{\end{facts}}
\newtheorem{definition}[thm]{Definition}
\newenvironment{definition*}{\begin{definition}\em}{\end{definition}}

\newcommand {\PP}{\mathbb{P}}
\newtheorem{exm}{Example}

\newcommand{\cE}{\mathcal{E}}
\newcommand{\cK}{\mathcal{K}}

\newcommand{\cF}{\mathcal{F}}

\newcommand{\cO}{\mathcal{O}}
\newcommand{\cQ}{\mathcal{Q}}
\newcommand{\cN}{\mathcal{N}}

\newcommand{\cI}{\mathcal{I}}

\DeclareMathOperator{\rk}{rk}

\DeclareMathOperator{\codim}{codim}

\DeclareMathOperator{\Pf}{Pf}

 \newcommand{\p}{\mathbf P}

\newcommand{\kk}{\mathbf{k}}

\SelectTips{cm}{} \newdir{ >}{{}*!/-5pt/@{>}}
\numberwithin{equation}{section}

\allowdisplaybreaks[1]
\begin{document}

\title{Vector spaces of skew--symmetric matrices of constant rank}

\author{Maria Lucia ~Fania} 
\email{fania@univaq.it}
\address{Dipartimento di Matematica Pura e Applicata, Universit\`{a} degli Studi dell'Aquila, \\Via Vetoio Loc. Coppito, 67100 L'Aquila, Italy}

\author{Emilia ~Mezzetti}
\email{mezzette@units.it}
\address{Dipartimento di Matematica e Informatica,
Universit\`{a} degli Studi di Trieste, 
Via Valerio 12/1, \\34127 Trieste, Italy}

\keywords{Skew--symmetric matrices, Grassmannian,  secant variety,  uniform bundle,  constant rank}
\subjclass[2000]{ 14M15,  15A30,  14N05,  14J60}

\thanks{The  authors were supported by MIUR funds, project \lq\lq\thinspace Geometria delle
variet\`a algebriche e dei loro spazi di moduli''. This work  was initiated during a visit of  the second author to MSRI with support of INdAM.} 

\begin{abstract}
We study the orbits of vector spaces of skew--symmetric matrices of constant rank $2r$ and type $(N+1)\times(N+1)$ under the natural action of $SL(N+1)$, over an algebraically closed field of characteristic zero. We give a complete description of the orbits for vector spaces of dimension $2$, relating them to some $1$-generic matrices of linear forms. We also show that, for each rank two vector bundle on $\PP^2$ defining a triple Veronese embedding of $\PP^2$ in $\GG(1,7)$, there exists a vector space of $8\times 8$ skew--symmetric matrices of constant rank $6$ whose kernel bundle is the dual of the given rank two vector bundle.
\end{abstract}





\maketitle

\section*{Introduction}

Vector spaces of skew--symmetric matrices of order $N+1$ and constant 
rank $2r$, with $N=2r+1$,  can be naturally interpreted as linear spaces contained
in the $(r-1)$-th secant variety  of the Grassmannian of lines $\GG(1,N)$, 
not meeting the $(r-2)$-th secant variety, which is its singular locus.  
Therefore   the special linear group  $SL(N+1)$ acts naturally on them and it 
is a natural problem to look for the maximal dimension of these spaces and to describe the orbits.

This problem has been considered both from the point of view of linear algebra and from that of algebraic geometry. An excellent survey of the results on the bounds, due to Ilic and Landsberg, is contained in 
\cite{IL}, in the wider context of matrices,  non--necessarily skew--symmetric ones. In particular, upper bounds are given for the dimensions of these spaces and a precise bound in a few cases. 

As for the orbits of these vector spaces not much is known so far. 
In \cite{MM} the case of $6\times 6$ skew--symmetric matrices of constant rank $4$ was considered and the orbits were completely classified, up to the action of $SL(6)$. The point of view adopted by  Manivel--Mezzetti is that of algebraic geometry:  to a vector space $M$ of matrices of constant rank 
$2r$ one can associate a vector bundle   
map 
\begin{eqnarray*}
\phi_M:  \cO_{\PP(M)}^{N+1}{\longrightarrow} \cO_{\PP(M)}(1)^{N+1}	
\end{eqnarray*}
on  the projective space  $\PP(M)$.  The kernel $\cK$ and   the image  $\cE$ 
 are vector bundles of ranks respectively $N+1-2r$ and $2r$, such that	$\cE$ is generated by its global sections,
	$\cE\simeq\cE^*(1)$ and
	 the splitting type of $\cE$ is $\cE\mid_l=
	\cO_l^r\oplus \cO_l^r(1)$, for all lines $l\subset M$. In  particolar $\cE$ is uniform and $c_{1}(\cE)=r=-c_{1}(\cK)$ (see Remark \ref{chern}.)

The classification in  \cite{MM} can be expressed  in terms of vector bundles on $\PP^1$ and $\PP^2$ with $c_1=2$.
Globally generated vector bundles on projective spaces with $c_1=2$ are completely described (see \cite{su}), in particular on $\PP^2$ there are $4$ of such rank two bundles and for each of them there is an orbit of skew--symmetric matrices of constant rank $4$ having it as dual of the  kernel bundle $\cK$.

In the present paper  we continue this study pursuing two objectives: on one hand the classification  
of vector spaces $M$ of dimension two, i.e. projective lines of  skew--symmetric matrices of any order, on the other hand that of  skew--symmetric matrices of rank $6$, that correspond to vector bundles with $c_1=3$.

To classify  the orbits  of projective lines of  skew--symmetric matrices of constant rank we rely on the fact that, in this case, congruence and strong equivalence of matrices are the same relation (see \cite{G}). This allows us to restrict our attention to the \lq\lq compression space'' matrices introduced by Eisenbud--Harris in \cite{EH}. The classification of the orbits,  given in Theorem \ref{classif}, is similar to the one of the rational normal scrolls, and follows from a link we establish between our matrices and $1$-generic matrices with two rows.

Globally generated vector bundles on projective spaces with $c_1=3$ which give a triple Veronese embedding of $\PP(M)$ in $\GG(1,N)$ have been studied in \cite{huh}. After refining such classification (see Theorem \ref{huh}), we  prove in Proposition \ref{quotientof} that the  non split vector bundles given in Theorem \ref{huh}  can all be expressed as  quotient of vector bundles of higher rank of a very particular form. This turns out to be crucial  in identifying some \lq\lq\thinspace building blocks'' skew--symmetric matrices, that we use to contruct  matrices of constant rank $6$ 
for each class of  rank $2$ bundles appearing in Theorem \ref{huh}.   This is done  by suitably projecting some direct sum matrices constructed  using the building blocks matrices of smaller rank (see Theorem \ref{proj}). We note that the case of rank $6$ is the first one in which infinitely many orbits appear.

As for $6\times 6$ matrices, also  in the case of $8\times 8$ matrices  there are no $\PP^3$'s of matrices of rank two less than the order. The first example of such a situation was given by Westwick (\cite{w2})  and it  is a $\PP^3$ of $10\times 10$ matrices of rank $8$. We  discuss this example in \S 5, where we also point out  the applications to the classification of degenerations of an important class of projective varieties, known as Palatini scrolls (see \cite{bazan-mezzetti}, \cite{faenzi-fania}).


\section{Definitions and preliminary results} \label{defs and prelim results}

Let $V$ be a vector space of dimension $N+1$ over the field $\kk$ 
(algebraically closed of characteristic $0$). We denote by $\GG(1,N)$ the Grassmannian 
of the vector subspaces of $V$ of dimension $2$, 
i.e. the projective subspaces of $\PP(V)$ of dimension $1$. $\GG(1,N)$ is embedded 
via the Pl\"ucker map in the projective space $\PP(\wedge^2V)$.
The group $SL(N+1)$, as well as  $PGL(N+1)$,
acts naturally on $\PP(V)$ and on $\PP(\wedge^2V)$.  If we fix a basis on $V$ then an element of $\wedge^2V$ can be thought as a skew--symmetric matrix and the action of $SL(N+1)$ on $\PP(\wedge^2V)$ is the congruence.
The orbits 
of the action on $\PP(\wedge^2V)$ are the 
Grassmannian and its secant varieties, and correspond respectively 
to the tensors of tensor rank $2, 4, \cdots, \left[\frac{N+1}{2}\right]$. 
$SL(N+1)$ (as well as  $PGL(N+1)$)
acts also naturally on the Grassmannian of the subspaces 
of $\PP(\wedge^2V)$ of any fixed dimension $d\le {{N+1}\choose {2}}$, 
or, in other words, on the skew--symmetric $(N+1) \times (N+1)$ matrices
of linear forms.
For this action, we are interested in describing the orbits of  
subspaces of constant rank, i.e. subspaces that are entirely contained 
in some orbit of the previous action on $\PP(\wedge^2V)$. Note that  the orbits of the action given by $SL(N+1)$ coincide with the orbits of the action given by  $PGL(N+1)$.

From now on, $M$ will denote a vector space of skew--symmetric matrices of order $N+1$, dimension $d$ and constant rank $\rk M=2r$.

If $N$ is odd and $2r$ is maximal, that is,  it is equal to $N+1$, then $d\le 1$, because the 
matrices of submaximal rank form a hypersurface 
of degree $\frac{N+1}{2}$, the Pfaffian, in $\PP(\wedge^2V)$. 
So we will assume that either $N$  is odd and
$2r$ is strictly less than $N+1$, or $N$ is even; in this last case the maximal rank is $N$ and the matrices of rank  $N-2$ have codimension $3$.  

\vskip 2mm
Given a vector space $M$ of matrices of constant rank 
$2r$ we can associate a vector bundle   
map 
\begin{eqnarray}\label{fasci}
\phi_M: V\otimes \cO_{\PP(M)}{\longrightarrow V}\otimes \cO_{\PP(M)}(1).	
\end{eqnarray}
on  the projective space    $\PP(M)$. 
 Since $M$ has constant  rank  then the kernel $\cK$,   the image  $\cE$ 
 and the cokernel  $\cN$ of  $\phi_M$ are vector bundles of ranks respectively $\rk \cK=\rk \cN=N+1-2r$, $\rk \cE=2r$ and determine short exact sequences: 
 \begin{eqnarray}
 \label{fibrativettassociati1}
0 \to \cK \to  V\otimes \cO_{\PP(M)}\to \cE \to 0	
\end{eqnarray}
\begin{eqnarray}
 \label{fibrativettassociati2}
0 \to \cE \to V\otimes \cO_{\PP(M)}(1) \to \cN \to 0	
\end{eqnarray}

\begin{prop} (\cite{EH}, \cite{IL})
\label{EHIL}
 \begin{enumerate}
	\item $\cE$ is generated by its global sections;
	\item  $\cE\simeq\cE^*(1)$;  $\cN\simeq\cK^*(1)$
	\item the splitting type of $\cE$ is $\cE\mid_l=
	\cO_l^r\oplus \cO_l^r(1)$, for all lines $l\subset M$. In  particolar $\cE$ is uniform.
	\end{enumerate}
\end{prop}
\begin{rem} \label{chern}
Note that  since the matrices are  skew--symmetric  the two short exact sequences (\ref{fibrativettassociati1})  and (\ref{fibrativettassociati2})   reduce to 
the single sequence 
 \begin{eqnarray}
 \label{fibrativettassociati}
0 \to \cK \to  V\otimes \cO_{\PP(M)}\to \cE \to 0	
\end{eqnarray}
with $\cE\simeq\cE^*(1)$. Moreover $c_{1}(\cE)=r=-c_{1}(\cK)$ and $\cK^{*}$ is generated by global sections.
\end{rem}

\begin{rem}
In the case $N=2r+1$, $M$ is contained in $S_{r-1}\GG(1,N)$, the top secant variety of
$\GG(1,N)$  strictly contained in $\PP(\wedge^2V)$. It
is naturally isomorphic to
${\overset\vee{\GG}}(1,N)$, the dual
of $\GG(1,N)$, which is the Pfaffian hypersurface.
Hence the Gauss map 
\begin{eqnarray}\label{mappaGauss}
\gamma: {\overset\vee{\GG}}(1,N) --\ra \GG(1,N)
\end{eqnarray}
is defined by the partial derivatives of the Pfaffian,
which are homogeneous polynomials of degree $r$.
\end{rem}

\begin{prop}
\label{gaussmap} Let $M$ be a vector space of dimension $d$ of $(N+1) \times (N+1)$ matrices of constant rank $2r$, with $N=2r+1$. Let $\gamma$ be the Gauss map in (\ref{mappaGauss}). Then $\gamma(\PP(M))$ is a Veronese variety  $v_r(\PP^{d-1})$ contained in the Grassmannian $\GG(1,N)$, or an isomorphic  projection of it.
\end{prop}

\proof
The  restriction of $\gamma$ to
$\PP(M)$ is regular, due to the hypothesis of constant rank, because $\PP(M)$ does not intersect  $S_{r-2}{{\GG}}(1,N)$, which is the indeterminacy locus of $\gamma$.
It remains to prove that $\gamma\mid_{\PP(M)}$ is biregular onto its image, i.e.  that $\PP(M)$ intersects a general fibre of $\gamma$ in only one point.  Let $l$ be a point of $\GG(1,N)$; the fibre $\gamma^{-1}(l)\subset {\overset\vee{\GG}}(1,N)$ can be interpreted as the set of hyperplanes containing the projective tangent space to the Grassmannian at $l$, $\TC_l$; hence $\gamma^{-1}(l)={\overset\vee\TC_l}$ is  a linear space and $\gamma^{-1}(l)\cap \PP(M)$ is also linear.

We can choose a basis $e_0, \cdots, e_N$ of $V$ such that $l=\langle e_0, e_1\rangle$; then the points in $\TC_l$ have Pl\"ucker coordinates $p_{ij}$, $0\leq i<j\leq N$ such that $p_{ij}=0$ for all $i\geq 2$.  Therefore ${\overset\vee\TC_l}$ is represented by matrices $a_{ij}$ whose first two rows and columns are zero; it can be seen as the linear span of a subgrassmannian $\GG(1,N-2)$, and the matrices of rank $2r$ are an open set in it whose complementar set is a hypersurface. If $\PP(M)$ intersects ${\overset\vee\TC_l}$ in positive dimension, then its intersection with this hypersurface is non-empty, so the rank is non-constant, a contradiction. Hence the map  $\gamma: \PP(M) \to \GG(1,N)$ is an embedding ($\gamma$ denotes also $\gamma_{|\PP(M)}$). 
\qed
\vspace{2mm}

\begin{facts}\label{facts}    
Let us recall some facts about embeddings in Grassmannians of lines (for details we refer to  \cite{artrento}).
Let $X$  be a smooth algebraic variety. To give 
a map $\varphi:X \to \GG(1,N)$ is equivalent to give a rank $2$ vector bundle $\cF$ on $X$  and an epimorphism $V\otimes \cO_{X} \to \cF\to 0$, where $V$ is an  $(N+1)$-dimensional subspace of $H^{0}(X,\cF)$. The map $\varphi$ is an embedding if any subscheme of $X$ of length two imposes at least three conditions to $V$.
Given an embedding in a Grassmannian $\varphi:X \to \GG(1,N)$  there is a ruled variety obtained by taking the union of all  lines defined by the points of $X$.  Consider $Y=\PP(\cF)$ the projective bundle associated to $\cF$,  let $\pi: Y\to X$ 
be the natural projection and let  $\cO_{\PP(\cF)}(1)$  be the  tautological line bundle on $Y$,  which has the property that $\pi_{*}(\cO_{\PP(\cF)}(1))=\cF$ and there is an epimorphism
 \begin{eqnarray*}\pi^{*}(\cF) \to \cO_{\PP(\cF)}(1)
  \end{eqnarray*}
 which induces an isomorphism $H^{0}(X,\cF)\cong  H^{0}(Y, \cO_{\PP(\cF)}(1))$.
 Hence an epimorphism  $V\otimes \cO_{X} \to \cF$ induces an epimorphism  $V\otimes \cO_{Y} \to \pi^{*}(\cF) \to \cO_{\PP(\cF)}(1).$  And viceversa, an epimorphism  $V\otimes \cO_{Y} \to  \cO_{\PP(\cF)}(1)$ induces an epimorphism  $V\otimes \cO_{X} \to \pi_{*}(\cO_{\PP(\cF)}(1))=\cF$. Thus a map $\varphi:X \to \GG(1,N)$  is equivalent to a map 
 $\overline{\varphi}:Y \to \PP^{N}$ of the corresponding ruled variety.  In general  $\overline{\varphi}$ is not an embedding. \\
 If the variety $X$ is two dimensional then $Y=\PP(\cF)$ is a $3$-fold and $\cO_{\PP(\cF)}(1)^{3}= c_{1}(\cF)^{2}-c_{2}(\cF)$ and hence we have 
 \begin{eqnarray} 
 \label{grado1Z}
 c_{1}(\cF)^{2}-c_{2}(\cF)=deg(\overline{\varphi})\cdot deg\ \overline{\varphi}(Y)
 \end{eqnarray}
Recall that it is also known the correspondence between globally generated vector bundle of rank $k$ and  maps to Grassmannians. In fact to give a regular map $\varphi_{V}: X \to \GG(k-1, N)$ is equivalent to give a  globally generated vector bundle $\cF$ of rank $k$  on $X$ and an  epimorphism $V\otimes \cO_{X} \to \cF\to 0$, where $V$ is an  $(N+1)$-dimensional subspace of $H^{0}(X,\cF)$.  Moreover, for any  epimorphism $V\otimes \cO_{X} \to \cF\to 0$ we get an exact sequence of vector bundles
\begin{eqnarray}
\label{duale}
0\to \cG \to V\otimes \cO_{X} \to \cF\to 0
\end{eqnarray}
hence  $V^{*}\otimes \cO_{X} \to \cG^{*}\to 0$ is also an epimorphism and thus  the pair $(\cG^{*},V^{*})$ defines a map $\phi_{V^{*}}: X \to \GG(N-k,N)$. The map $\phi_{V^{*}}$ is said the dual map of $\varphi_{V}$. It follows that $\varphi_{V}(X)$ and $\phi_{V^{*}}(X)$ are naturally isomorphic by the duality $\GG(k-1,N)\cong \GG(N-k, N)$, see \cite{su} for details. 
\end{facts}

\begin{prop}
\label{-K da gaussmap} 
 If $M$ is a  vector space of $(N+1) \times (N+1)$ matrices of constant rank $2r$, with  $2r=N-1$, then $\cK^{*}$ is a vector bundle of rank $2$ which defines an $r$-tuple embedding of   $\PP(M)$ in $\GG(1,N)$.
\end{prop}
\proof
It is enough to note that the Gauss map (\ref{mappaGauss}) restricted to $\PP(M)$ is given by the rank $2$ bundle $\cK^{*}$ on $\PP(M)$ and the epimorphism $V\otimes \cO_{\PP(M)} \to \cK^{*}$ is  obtained by dualizing   (\ref{fibrativettassociati}). 
\qed


\section{Vector spaces of skew--symmetric matrices of dimension two} \label{rette}

In this section we  study the first non-trivial case of vector space of skew--symmetric matrices of constant rank $2r$ and of order $N+1$, that is the case in which  $\dim M=2$. We will give a complete classification of the orbits and an explicit description of the corresponding matrices. This will be possible because the vector bundles on $\PP^1$ are all decomposable.

We recall the notion of strictly equivalent pencils of matrices. 
\begin{definition} Two pencils of matrices $aA+bB$  and  $aA_{1}+bB_{1}$  are called strictly equivalent if there exist two non singular matrices $P, Q$ with entries in $\kk$ such that 
 \begin{eqnarray*}P(aA+bB)Q=aA_{1}+bB_{1} \end{eqnarray*}
\end{definition}
The following well known Theorem says, in particular, that for pencils of skew--symmetric matrices the notion of ``strictly equivalent"  coincides with the notion of ``congruent".
\begin{thm}  (\cite{G}, Chap. XII, Theorem 6) \label{Gantmacher} Two strictly equivalent pencils of complex symmetric (or skew--symmetric) matrices are always congruent. 
\end{thm}
This theorem allows to use, in the case of $\dim M=2$, the result obtained by Eisenbud-Harris in \cite{EH}, where they consider the classification of vector spaces of matrices of linear forms for the relation of strict equivalence. We recall  some terminology from \cite{EH}. 
\begin{definition} Let $M$ be a vector space of matrices of constant rank. $M$  is {\it nondegenerate}  if the kernels 
of the matrices of $M$ intersect in the 
zero subspace and the images of the elements of $M$
generate the vector space $V$. \end{definition}

This is equivalent to say that  $M$ is not 
$SL(N+1)$-equivalent to a space of matrices with a row 
or a column of zeroes. In other words $N$ is 
the minimum integer such that $M$ can be embedded in 
$\langle \GG(1,N)\rangle$. 

From now on we 
will consider only nondegenerate vector spaces of $(N+1) \times (N+1)$ matrices of constant rank $2r$. 

\begin{definition} Let $M$ be a vector space of matrices of constant rank $2r$.  $M$ is a {\it compression space} if there exist subspaces 
$V', W'\subset V$, 
such that every matrix in $M$ maps $V'$ into $W'$  and $\rk M = 2r = \codim V'+\dim W'$. 
\end{definition}

It is easy to see, by an appropriate choice of basis of $V$, 
that $M$ is a compression space if and only if  it is 
$SL(N+1)$-equivalent to a space of $(N+1)\times(N+1)$ matrices having 
a common block of zeroes of size $(N+1-h)\times(N+1-k)$ with $h+k=2r$.

\begin{prop} (\cite{EH}, Corollary 2.2) \label{compression}
If $\dim M=2$ then $M$ is a compression space.
\end{prop}

The kernel bundle $\cK$ of the map $\phi_M$  (see (\ref{fasci})) is of the following form
\begin{equation}
	\cK=\cO_{\PP^1}^{m_0}\oplus \cO_{\PP^1}(-1)^{m_1}
	\oplus \cdots \oplus \cO_{\PP^1}(-k)^{m_k}
\end{equation}
where $m_0, m_1, \cdots, m_k$ are non-negative integers. 

\begin{prop}\label{bound_N}
  If $M$ is nondegenerate, $rk M=2r$ and   $\dim M=2$, then $2r\leq N \leq 3r-1$.
  \end{prop} 
  
  \proof
 Comparing the ranks and the first Chern 
 classes of the bundles appearing in (\ref{fibrativettassociati}), 
 we get: 
$m_0+m_1+\cdots +m_k=N+1-2r$ and $m_1+2m_2+
\cdots +km_k=r$.
The assumption that $M$ is nondegenerate 
implies moreover that $m_0=0$. The thesis follows by 
computing the dimensions of the cohomology 
groups of (\ref{fibrativettassociati}), taking into account that 
$\rk \cK\geq 1$.
  \qed
  \vspace{2mm}
  
  The splitting type of $\cK^*$  is a partition of $r$ of the form $r=r_1+r_2+\cdots +r_h$. 
  From Proposition \ref{bound_N}, it follows that  
  the length $h$ of the partition is $h=N+1-2r$.
    
  Conversely, for each $r\geq 1$ and each partition of $r$, we shall exhibit 
  a pencil of matrices  having it as splitting type of the associated bundle $\cK^{*}$. 
  All the 
  corresponding matrices will have a $(N+1-r)\times (N+1-r)$ block of zeroes.
  \vspace{3mm}
  
  Let us start  with a few examples.  We fix a basis $e_0, \cdots e_N$  for  the vector space $V$ and let $M$ be  a pencil of skew--symmetric matrices. Then a general matrix  in $M$ is of the form $aA+bB$, with $A, B$ skew--symmetric matrices with constant entries and $a,b\in \kk$.

  \begin{egs}\label{esempi}$\\$ We will always write our skew--symmetric matrices indicating only the entries  on the strict upper triangular part.
  
	\quad $\bullet \ r=1$
	
	In this case $\rk M=2$, there is only one orbit, $M$ is contained in 
	$\GG(1,2)$ and corresponds to a pencil of lines in the plane. 
	An element of the orbit
	is  the following $3\times 3$ matrix:
	
  \begin{eqnarray*}\begin{pmatrix}
          a & b  \\
          & 0  \\
       \end{pmatrix} \end{eqnarray*}
  \\   
  $\bullet \ r=2$ 
  
   In this case $\rk M=4$, there are two orbits, corresponding to the 
   bundles $\cO_{\PP^1}(2)$ and $\cO_{\PP^1}(1)\oplus \cO_{\PP^1}(1)$, 
   formed by matrices of order $5\times 5$, resp. $6\times 6$ (see \cite{MM}, \S 3):
  \begin{eqnarray*}\begin{pmatrix}
           0 & a & b  & 0  \\
           & 0 & a  & b  \\
         &   & 0  & 0  \\
          &  &    & 0  \\
       \end{pmatrix}, \  \  \   \begin{pmatrix}
         0 & a & b & 0 & 0  \\
        & 0 & 0 & a & b  \\
          &  & 0 & 0&0  \\
       & &  & 0 & 0  \\
        &  & &  & 0  \\
       \end{pmatrix} \end{eqnarray*}
    \\   
  $\bullet \ r=3$ 
  
In this case $\rk M=6$.
We give three examples which we denote by $M_7, M_8, M_9$, respectively. They correspond to the 
bundles $\cO_{\PP^1}(3)$, $\cO_{\PP^1}(2)\oplus\cO_{\PP^1}(1)$, $\cO_{\PP^1}(1)\oplus\cO_{\PP^1}(1)\oplus\cO_{\PP^1}(1)$; 
the orders of the matrices are  $7\times 7$, $8\times 8$ and $9\times 9$.

 \begin{eqnarray*}
 M_7\!=\!\begin{pmatrix}
          0 & 0 & a & b & 0 & 0 \\
         & 0 & 0 & a & b & 0 \\
         &  & 0 & 0 & a & b \\
        & &  & 0 & 0 & 0 \\
  & &  &  & 0 & 0 \\
         & & & &  & 0 \\
   \end{pmatrix}\!\!, 
   M_8\!=\!\begin{pmatrix}
          0 & 0 & a & b & 0 & 0 & 0\\
         & 0 & 0 & a & b & 0 & 0\\
          &  & 0 & 0 & 0 & a & b \\
   & &  & 0 & 0 & 0 & 0\\
&  & &  & 0 & 0 & 0\\
         & & & &  & 0 & 0\\
         &  &  & & &  & 0\\
       \end{pmatrix}\!\!, M_9\!=\!\begin{pmatrix}
        0 & 0 & a & b & 0 & 0 & 0 & 0\\
        & 0 & 0 & 0 & a & b & 0 & 0\\
         & & 0 & 0 & 0 & 0 & a & b \\
          &  &  & 0 & 0 & 0 & 0 & 0\\
        &&  &  & 0 & 0 & 0 & 0\\
       & & & &  & 0 & 0 & 0\\
       & & && &  & 0 & 0\\
       & & & &  & &  & 0\\
       \end{pmatrix}
        \end{eqnarray*}

  \end{egs}
  
  \begin{thm}\label{gli_esempi}
  Let $r=r_1+r_2+\cdots +r_h$ be a partition of $r$, $r_1\geq r_2\geq\cdots \geq r_h$. 
  There exists an orbit of pencils of 
  skew--symmetric matrices of constant rank $2r$ and order $N+1$ with  $N=2r+h-1$, whose associated bundle
  $\cK^{*}$ has splitting type $(r_1, \cdots, r_h)$.
   \end{thm}
   \proof
   In the orbit there is a matrix $F$ of 
   the following type. 
    \begin{equation}\label{matrixF}
F=\begin{pmatrix}
        0_{r} &\overline{F} \\
        - \overline{F}^t &  0_{N+1-r} \\
       \end{pmatrix}
       \end{equation}
   where $0_k$ denotes the zero matrix of order $k$ and $\overline{F}$ is a block matrix of the form 
   \begin{equation}
   \label{Fbar}
 \overline{F}=\begin{pmatrix}
         U_{r_1}&  &  &  &  \textrm{\large{0}}  \\
                      &U_{r_2}&&&\\
         & & \ddots & & \\
         0& && & U_{r_h} \\
       \end{pmatrix}.
       \end{equation}
   For every $i=1, \ldots, h$, $U_{r_i}$ is of type $r_i\times(r_i+1)$ and
   \begin{equation}\label{blocchi}
U_{r_{i}}=\begin{pmatrix}
         a & b & 0 & 0 & \cdots & . & .&0\\
         0 & a & b & 0 & \cdots & . &.&0\\
         0 & 0 & a & b & \cdots & . &.&0\\
          &  &  & & \ddots & &  &  \\
         . & . & . & . & .& a & b & 0 \\
         . & . & . & . & .& . & a & b \\
       \end{pmatrix}
       \end{equation}
       
   The resulting matrix $F$ has clearly 
   constant rank $2r$ and the associated bundle
   $\cK^{*}$ is as required. 
  
   \qed
   \vspace{2mm}

   We  recall  few facts about $1$-generic matrices  which  will be used in classifying pencils of skew--symmetric matrices of constant rank. We refer to \cite{H} for the definition and 
 properties of $1$-generic matrices.

  \begin{definition} Let $\Omega$ be a matrix of linear forms on  $\PP^{n}$. We say that $\Omega$ is {\it $1$-generic}  if no matrix of linear forms conjugate to $\Omega$ has a zero entry.
  \end{definition}
  
 \begin{prop}(\cite[Prop. 9.12 and its generalization]{H}) 
 \label{harris}
 \begin{enumerate}
 \item Any $1$-generic $2\times (n-1)$ matrix $\Omega$ of linear forms on $\PP^{n}$   is conjugate for some $\ell$ to the matrix 
\begin{eqnarray*}\Omega_{0}=\left( \begin{array}{cccccccccc}
  z_{0}& \cdots & z_{\ell -1} & z_{\ell +1}&  \cdots &z_{n -1}\\
 z_{1}& \cdots & z_{\ell } & z_{\ell +2}&  \cdots &z_{n }
 \end{array}
\right)
\end{eqnarray*}
where $z_0,\cdots, z_n$ are  homogeneous coordinates  on $\PP^{n}$. 
\item Let $\Omega$ be a  $1$-generic $2\times k$  matrix  of linear forms on $\PP^{n}$,  $k\leq n-1$,
whose entries   span $V^*$. Then for some sequence of integers $a_1, \cdots, a_{\ell}\ (\ell=n-k), \Omega$ is conjugate to the matrix
 \begin{eqnarray*}\label{catalecticantblocks}
\Omega_{a}=\left( \begin{array}{cccccccccccccccccc}
  z_{0}& \cdots & z_{a_{1} -1} &\vline& z_{a_{1} +1}&  \cdots &z_{a_{2} -1} &\vline& z_{a_{2} +1}&  \cdots &z_{a_{\ell} -1} &\vline& z_{a_{\ell} +1}&  \cdots &z_{a_{n} -1} \\
 z_{1}& \cdots & z_{a_{1}  } &\vline& z_{a_{1}+2}&  \cdots &z_{a_{2}} &\vline& z_{a_{2} +2}&  \cdots &z_{a_{\ell} } &\vline&z_{a_{\ell} +2}&  \cdots &z_{a_{n}} 
 \end{array}
\right) ,
 \end{eqnarray*}
that is to a matrix consisting of $\ell +1$ blocks of size $2\times a_1, \cdots, 2\times (n-a_{\ell}-1)$ with each block a {\it catalecticant}, that is a matrix in which  $a_{i,j+1}=a_{i+1,j}$ for all $i, j$. 
\end{enumerate}
\end{prop}
   \vspace{2mm}

\begin{rem} Computing the kernel of the 
matrix $F$ in  (\ref{matrixF}), one gets the family of 
$(h-1)$-spaces of a rational normal scroll 
of type $(r_1, \cdots, r_h)$. 

\noindent In fact  if we write the vector $X=(x_0,\cdots, x_{r-1},x_{r},\cdots, x_{N})$ as $X=(X^{r},\overline{X})$, where
 $X^{r}=(x_0,\cdots, x_{r-1}),\,  \overline{X}=(x_{r},\cdots, x_{N}),$  then 
 \begin{eqnarray*}
 ker(F)&=&\{X=(x_0,\cdots, x_{r-1},x_{r},\cdots, x_{N}) \, | \, F\cdot X=0\}\\
 &=&\{X=(X^{r},\overline{X}) \, | \, \overline{F}\cdot \overline{X}=0\quad  \text{and} \quad \overline{F}^t \cdot X^{r}=0\}
 \end{eqnarray*}
 From $\overline{F}^t\cdot X^{r}=0$ we get that $x_0=\cdots=x_{r-1}=0$. 
 
Moreover,  from (\ref{Fbar}) 
    we deduce that $\overline{F}\cdot \overline{X}=0$
       is equivalent to 
       \begin{equation}\label{leU} U_{r_{1}}X_{r_{1}}=0, \cdots,  U_{r_{h}}X_{r_{h}}=0,\end{equation}
        where 
  $X_{r_1}=(x_{r}, \cdots, x_{r+r_{1}})$, 
  $\cdots$, $X_{r_{h}}=(x_{r+r_{1}+r_{2}+\cdots +r_{h-1}+h-1}, \cdots, x_{r+r_{1}+r_{2}+\cdots +r_{h}+h-1})=$
  $(x_{r+r_{1}+r_{2}+\cdots +r_{h-1}+h-1}, \cdots, x_{N})$. Hence  the conditions defining $ker(F)$ 
are equivalent to the fact that  the following matrix has rank $1$:
 \begin{eqnarray*}\label{catalecticant}
\Omega(\overline{X})=\left( \begin{array}{cccccccccc}
 x_{r} &x_{r+1}& \cdots & x_{r+r_{1}-1} & \vline & x_{r+r_{1}+1} & x_{r+r_{1}+2}& \cdots & \vline & \cdots\\
  x_{r+1} &x_{r+2}& \cdots & x_{r+r_{1}} & \vline & x_{r+r_{1}+2}&x_{r+r_{1}+3} &\cdots & \vline& \cdots
\end{array}
\right). 
  \end{eqnarray*}
 The matrix $\Omega(\overline{X})$ consists of $h$ blocks of size $2\times r_{1}, \cdots, 2\times r_{h},$ respectively, with each block catalecticant. If  we denote such blocks by $\Omega_i(X)$, with $i=1, \cdots, h$, each block gives a rational normal scroll, that is  the determinantal variety defined by $\text{rank}\ \Omega_i(X)=1$. 
 \end{rem}

 \begin{thm}\label{classif}
 Let $M$ be a nondegenerate vector space of 
 dimension $2$ of matrices 
 of constant rank $2r$ and order $N+1$. 
 Then $M$ is $SL(N+1)$-equivalent by congruence and strict equivalence
 to one of the matrices  of Theorem \ref {gli_esempi}.
 \end{thm}
 
 \proof
 If we think of $M$ 
as a subspace of $\wedge^2V$,  we can write it as  $M=\left\langle \omega, \omega'\right\rangle$, 
 where $\omega, \omega'$ are tensors of tensor rank 
 $r$ whose linear combinations have all rank equal to $r$.
 There exist expressions 
 \begin{equation}\label{omega}
\omega=u_0\wedge v_0+\cdots+u_{r-1}\wedge v_{r-1}, \  
 \omega'=z_0\wedge w_0+\cdots +z_{r-1}\wedge w_{r-1}.\end{equation}
 Let $L$ and  $L'$ be the
subspaces of $V$ generated by the vectors
$u_0, v_0, \cdots, u_{r-1}, v_{r-1}$ and
$z_0, w_0, \cdots, z_{r-1}, w_{r-1}$, respectively. Note that $\dim
L=\dim L'=2r$, because
$\omega\in\wedge^2L$,  $\omega'\in\wedge^2L'$
 and $2r$ is the minimal dimension of a vector space such that there exist
 skew--symmetric tensors of tensor rank $r$. 
 So the given generators are linearly independent.
 
 Since $M$ is a compression space by Proposition 
 \ref{compression}, 
 there exist vector subspaces $V', W'$ of $V$ 
 such that  every matrix in $M$ maps $V'$ into $W'$  and $\rk M = 2r = \codim V'+\dim W'$. 
 
 The strategy of the proof is to first analyze 
 the case $\codim V'=\dim W'=r$: we will  
 show that $M$ is in the orbit of one of the 
 examples of Theorem \ref{gli_esempi}. 
 Then we will consider the other possibilities for 
 $\codim V'$ and $\dim W'$ and will prove that they cannot
 occur.
 
 $\bullet$ \ Assume $\codim V'=\dim W'=r$. $\\$
 We can choose a basis of $V$ such that $V'=\left\langle e_r,
 \cdots ,e_N\right\rangle$ and $W'=\left\langle e_0,
 \cdots ,e_{r-1}\right\rangle$. This means that 
 the submatrix of the last $N+1-r$ rows and columns
 is the zero matrix.
 Therefore in the expression (\ref{omega}) we have 
  $u_0,\cdots,u_{r-1}, z_0, \cdots,z_{r-1}\in W'$. Possibly changing basis
  in $W'$, we can assume that 
  \begin{equation}\label{riomega}
	\omega=e_0\wedge v_0+\cdots+e_{r-1}\wedge v_{r-1}.\end{equation} 
		Therefore $z_j=\sum_{i=0}^{r-1}\lambda_{ji}e_i$, 
  for all $j$ and suitable scalars $\lambda_{ji}$. 
  Using bilinearity, we can assume that 
  $\omega'$ has the form 
 \begin{equation}\label{omega'}
 \omega'=e_0\wedge w_0+\cdots+e_{r-1}\wedge w_{r-1}.\end{equation} 
 
 Let  $\overline{M}$ be a matrix in $M$, it can be written as $aQ+bP$, where $Q$ and $P$ are the matrices representing $\omega$ and $\omega'$ respectively. Let us denote its general element
 by $aq_{ij}+bp_{ij}$. Moreover $\overline{M}$ has the following form:
 \begin{equation}\label{tilde}
\begin{pmatrix}
         M' & \tilde{M} \\
         -\tilde{M}^t &  0 \\
       \end{pmatrix}
       \end{equation}
 
\noindent  where $M'$ and $\tilde{M}$ are matrices of linear forms in  $a,b$
 of size respectively $r\times r$ and $r\times (r+h)$.  
 Note that $\tilde{M}$ has
  maximal rank $r$ for every pair $(a,b)$.  
  For each $(a,b)$, 
  the vectors of   
   $\ker \overline{M}$ 
 have the form $(0,\cdots,0,x_r,\cdots, x_N)$, where
 $(x_r,\cdots, x_N)$ belongs to the kernel of $\tilde{M}$, 
 which is a vector space 
 of dimension $N-2r+1=h$. Letting $(a,b)$ vary, we obtain a variety $Y$ of dimension $h$ in $\PP^{N-r}$ covered by linear spaces (see Facts \ref{facts}). 
 The equations of $Y$ are the $2\times 2$ minors 
 of the $2\times r$ matrix 
 \begin{equation}\label{1-generica} \Pi=
\begin{pmatrix}
          \tilde{Q}_1X'& \cdots & \tilde{Q}_rX' \\
         \tilde{P}_1X' &  \cdots & \tilde{P}_rX' \\
       \end{pmatrix}
       \end{equation}
 \noindent where $\tilde{Q}$ and $\tilde{P}$ are for $\tilde{M}$ the analogous of $Q$ and $P$  for $\overline{M}$, and $X'$ is the column matrix with entries 
 $x_r,\cdots, x_N$; moreover $\tilde{Q_1}, \cdots, \tilde{Q_r}$ are the rows of 
 $\tilde{Q}$ and similarly for $\tilde{P}$.
 The matrix $\Pi$ is $1$-generic because $\tilde{M}$ has 
 constant rank $r$ and its entries 
 generate $(\kk^{N-r-1})^* $ because 
 $\overline{M}$ is nondegenerate. Hence $Y$ is a rational normal scroll in $\PP^{N-r}$.
 
	 By Proposition \ref{harris} it follows that $\Pi$ is conjugate 
 to a matrix $\Pi'$ with $h$ blocks:
 \begin{equation}\label{blocchi2}
 \Pi'=\begin{pmatrix}
          \Pi_1& \Pi_2 & \cdots & \Pi_h 
       \end{pmatrix}
       \end{equation}
 \noindent  where each $\Pi_i$ is a catalecticant matrix. This means that 
 $\Pi'$ is obtained from $\Pi$ by suitably multiplying it at the left and at the right
 by invertible scalar matrices. 
 
 By direct computations, one checks that left 
 multiplication of $\Pi$ by a $2\times 2$ matrix corresponds to changing
 generators for the pencil $aQ+bP$, and right 
 multiplication by a $r\times r$ matrix corresponds to 
 changing the last $N-r+1$ vectors of the basis of $V$.
 This shows that the matrix $\overline{M}$ is equivalent to a matrix of the desired form.

$\bullet$ \ Assume $\codim V'=r-k, \dim W'=r+k$, $k>0$. 

We choose a basis $(e_0, \cdots, e_N)$ of $V$ such that $V'=\left\langle e_{r-k},\cdots, 
e_N\right\rangle$ and $W'=\left\langle e_0,\cdots, e_{r+k-1}\right\rangle.$ 
In view of  skew--symmetry, the matrix $\overline{M}$ is concentrated 
in the first $r-k$ rows and columns, except for a 
$2k\times 2k$ skew--symmetric submatrix $D$ in the rows and 
columns of indices $r-k,\cdots, r+k-1$, as sketched in 
(\ref{compression_k}): 
\begin{equation}\label{compression_k}
 \overline{M}=\begin{pmatrix}
          A    & B & C \\
          -B^t & D & 0 \\
          -C^t & 0 & 0 \\
       \end{pmatrix}
       \end{equation}

\noindent Note that, since $\rk \overline{M}=2r$, necessarily $\rk C=r-k$ and $\rk D=2k$. 
But $D$ is a skew--symmetric matrix of order $2k$ whose 
entries are linear forms in $a,b$, hence its Pfaffian  vanishes for some pair $(a,b)\neq (0,0)$ (because the base field $\kk$ is algebraically closed).
This contradicts the assumption that the rank of $\overline{M}$ is constant and equal to $2r$.
 \qed
 \vspace{2mm}
 
   \begin{corol} The orbits of vector spaces of dimension 
 two of matrices of constant rank $2r$ and
 order $N+1$ are the ones of Theorem \ref{classif} 
and those of nondegenerate spaces 
 of lower order with a suitable number of zero rows and columns added.
 \end{corol}
 
 \subsection{An algorithm to compute the dimension of the orbits} \label{algo}
 It would be interesting to compute the dimension of the orbits. In the case of lines of skew--symmetric matrices of rank $4$, that is $r=2$, this has been done in \cite{MM}. 
 
 To compute the dimension for $r\geq 3$, we use the computer algebra system
 Macaulay2  (\cite{macaulay}) with the script suggested to us by Giorgio Ottaviani, and we thank him for this.

We consider here the examples of the orbits with $r=3$.
We write down the case $M_7$.  For ${M_7}', {M_7}'', M_8,  {M_8}', M_9$ one makes the appropriate changes.  For $i=7,8$, with  ${M_i}'$ we denote the matrix obtained from $M_i$ by adding   one row  and one column (the last ones) of zeroes and with  ${M_i}''$ we denote the matrix obtained from $M_i$ by adding  two rows and two columns (the last ones) of zeroes.
 \vspace{2mm}

{\textrm
R=QQ[a,b]
 
--M is our matrix 

N=map(R$^7$,R$^7,\{(0,3)=>a,(1,4)=>a,(2,5)=>a,(0,4)=>b,(1,5)=>b,(2,6)=>b\})$

M=N-transpose(N)

$P=(M)^{\{0\}}_{\{1..6\}}$

for s from 1 to 5 do P=P$|$(M)$^{\{s\}}_{\{(s+1)..6\}}$

-- we create  P with  21 components which represents the matrix  M

for i from $0$ to $6$ do

for j from $0$ to $6$ do

E$_{(i,j)}$=map(R$^7$,R$^7,\{(i,j)=>1_{R}\})$

-- E$_{(i,j)}$ are the elementary matrices

W=(transpose(E$_{(0,0)}$)*M+M*E$_{(0,0)})^{\{0\}}_{\{1..6\}}$

for s from 1 to 5 do W=W$|$(transpose(E$_{(0,0)}$)*M+M*E$_{(0,0)})^{\{s\}}_{\{(s+1)..6\}}$

WW1=sub(P,$\{a=>1_R,b=>0_R\}$)$||$sub(W,$\{a=>0_R,b=>1_{R}\})$

WW2=sub(W,$\{a=>1_R,b=>0_R\}$)$||$sub(P,$\{a=>0_R,b=>1_{R}\}$)

Z=exteriorPower(2,WW1)+exteriorPower(2,WW2)

-- now  Z has 210 components and represents the derivative of the action of E$_{(0,0)}$

-- with the following commands we repeat the above for $49$ times obtaining a matrix $50\times210$, 

where the first two rows are equal

for i from $0$ to $6$ do for j from $0$ to $6$ do

$\{$W=(transpose(E$_{(i,j)}$)*M+M*E$_{(i,j)}$)$^{\{0\}}_{\{1..6\}},$

for s from $1$ to $5$ do W=W$|$(transpose(E$_{(i,j)}$)*M+M*E$_{(i,j)})^{\{s\}}_{\{(s+1)..6\}},$

WW1=sub(P,$\{a=>1_R,b=>0_{R}\}$)$||$sub(W,$\{a=>0_R,b=>1_{R}\}),$

WW2=sub(W,$\{a=>1_R,b=>0_{R}\}$)$||$sub(P,$\{a=>0_{R},b=>1_{R}\}$),

Z=Z$||$exteriorPower(2,WW1)+exteriorPower(2,WW2)
$\}$

Z;

rank(Z)

--  rank(Z) represents the affine dimension of the orbit.

 In this case we get rank(Z)=$39$.}

  \vspace{5mm}

Recall that $N=\dim  \PP(V)$.
 
 If $N=6$ there is only one orbit. We get that the orbit  $\cO_{7}$ of $M_{7}$ has  $\dim {\cO}_{7}=38$. Hence it is open in $\GG(1,\PP(\wedge^2(\CC^{7}))$ and its complementary is formed by the lines which intersect $S_{1}(\GG(1,6)).$
 
  If $N=7$ there are two orbits. One,  ${\cO}_{7}'$, corresponding to  
  ${M}_{7}'$. The other,  ${\cO}_{8}$, corresponding to $M_8$.
  
We get that $\dim {\cO}_{8}=47$.  This is the expected dimension. Indeed $\dim\GG(1,\PP(\wedge^2(\CC^{8}))=52$ and $\deg S_2\GG(1,7)=4$,  because its equation is the Pfaffian of a $8\times 8$ matrix of linear forms. As for the dimension of  ${\cO}_{7}'$, we have that 
  $\dim {\cO}_{7}'=  \dim {\cO}_{7}+ \dim \check{\PP}^7=38+7=45,$ because a matrix of   ${\cO}_{7}'$ determines in a unique way a hyperplane in ${\PP}^7.$
  
   If $N=8$ there are three orbits,  ${\cO}_{7}''$, ${\cO}_{8}'$, ${\cO}_{9}$, whose dimensions are respectively $\dim {\cO}_{7}''=  \dim {\cO}_{7}+ \dim {\GG}(6,8)=38+14=52$, because a matrix of   ${\cO}_{7}''$ determines in a unique way a codimension two subspace of ${\PP}^8,$ $ \dim {\cO}_{8}'=\dim {\cO}_{8}+ \dim \check{\PP}^8=47+8=55,$     $\dim {\cO}_{9}=56.$ 

   \vspace{5mm}

\section{Building blocks matrices} \label{dim maggiore}

We turn now to vector spaces $M$ of skew--symmetric matrices of dimension at least $3$, where the situation is much more complex. We will consider therefore mainly the cases of low rank, and precisely those of rank $2r\leq6$, because the vector bundles on the projective spaces which are globally generated are classified for $c_1\leq2$ (\cite{su}) as well as for $c_1=3$ and rank $2$ (\cite{huh}).

\subsection{Rank $2$}
For $r=1$, we get the classification of the linear spaces contained in a Grassmannian of lines $\GG(1,N)$. It is well known that the maximal ones belong to one of the following two types:
\begin{itemize}
\item[(i)] the lines contained in a fixed $\PP^2$;
\item[(ii)] the lines passing through a fixed point in $\PP^N.$
\end{itemize}

In case (i) the corresponding exact sequence of bundles is
\begin{equation}0\to\cO_{\PP^{2}}(-1)\to \cO_{\PP^{2}}^3\to T_{\PP^{2}}(-1)\to 0\end{equation}
and a matrix in the orbit is
\begin{equation} \label{primo}\begin{pmatrix}
          a & b \\
          & c\\
       \end{pmatrix}.\end{equation}
In case (ii), we get a $\PP^{N-1}\subset\GG(1,N)$, the exact sequence is
\begin{equation} 0\to\Omega_{\PP^{N-1}}(1)\to \cO_{\PP^{N-1}}^{N+1}\to \cO_{\PP^{N-1}}\oplus\cO_{\PP^{N-1}}(1)\to 0\end{equation}
and a representative matrix is
\begin{equation}\label{secondo}\begin{pmatrix}
     a_1 & \dots & a_N\\
      &\dots & 0\\
       && \vdots\\
     && 0\\
       \end{pmatrix}.\end{equation}

\subsection{Rank $4$}
  For $r=2$, a classification of the orbits for matrices of order at most $6\times 6$ has been given in \cite{MM}. The result is that there are no vector spaces of dimension $4$ of such matrices (see also \cite{w2}), while the orbits of vector spaces of dimension $3$ are completely described.   In the case of $5\times 5$ matrices, there is only one orbit, i.e. the open subset of $\GG(2,9)$ complementar  to the irreducible subvariety of codimension $1$ representing $2$-planes meeting $\GG(1,4)$. The exact sequence (\ref{fibrativettassociati1}) in this case is
\begin{equation}\label{elencw}  0\to\cO_{\PP^{2}}(-2)\to \cO_{\PP^{2}}^5\to \cE\to 0
\end{equation} 
where $\cE$ is an indecomposable uniform bundle of rank $4$. A representative matrix in $M$ is 
\begin{equation}\label{O2} \begin{pmatrix}
        0&0 &a & b \\
      &a & b & c \\
  &&c & 0 \\
&  & &0 \\
         \end{pmatrix}.
           \end{equation} 
As for $6\times 6$ matrices,    there are $4$ orbits, one for each   of the globally generated rank 
$2$ bundles on $\PP^2$ with $c_1=2$ that are: $\cO_{\PP^{2}}\oplus\cO_{\PP^{2}}(2)$, $\cO_{\PP^{2}}(1)\oplus\cO_{\PP^{2}}(1)$, the restricted null-correlation bundle, which is a quotient of $\cO_{\PP^{2}}(1)\oplus T_{\PP^2}(-1)$, and the Steiner bundle, which is a quotient of $T_{\PP^2}(-1)\oplus T_{\PP^2}(-1)$.
The corresponding image bundles $\cE$ are respectively:  $\cE$ appearing in \ref{elencw}, $T_{\PP^2}(-1)\oplus T_{\PP^2}(-1)$, $\cO_{\PP^{2}}\oplus\cO_{\PP^{2}}(1)\oplus T_{\PP^2}(-1)$ and $\cO_{\PP^{2}}^2\oplus\cO_{\PP^{2}}(1)^2$. Representative matrices in the last $3$ cases are for instance:
         \begin{equation}  
     \begin{pmatrix}
      a&b &0 & 0 & 0\\
       &c & 0 & 0 & 0\\
      &&0 & 0 &0\\
      &  &&a&b \\
        & & &&c\\
         \end{pmatrix}, 
     \begin{pmatrix}
   0&0 &a & b & c\\
     &a & b & 0 & 0\\
&&c & 0 &0\\
&& &0&0 \\ 
  &&&&0\\
     \end{pmatrix},
     \begin{pmatrix}
       0&0 &a & b & c\\
      &a&b&c&0\\
     &&0&0&0\\
     &&&0&0\\
    &&&&0\\
      \end{pmatrix}.  
\end{equation}
For matrices of order at least $7$, the possible image bundles $\cE$ remain the same, whereas the dual of the kernel, $\cK^{*}$,  can be either $\cO_{\PP^{2}}(1)\oplus T_{\PP^2}(-1)$ or $T_{\PP^2}(-1)\oplus T_{\PP^2}(-1)$ or a quotient of it of rank $3$ (up to trivial direct summands). Examples of matrices are the following:
 \begin{equation}  
     \begin{pmatrix}
 0&0 &0& a & b & c\\
     &a & b & 0 & 0 &0\\
&&c & 0 &0&0\\
  &   & &0&0 &0\\ 
   &&&&0&0\\
     &&&&&0\\
     \end{pmatrix},
     \begin{pmatrix}
   0&0 &0& 0&a & b & c\\
    &a&b&c&0&0&0\\
    &&0&0&0&0&0\\
     &&&0&0&0&0\\
    &&&&0&0&0\\
    &&&&&0&0\\
 &&&&&&0\\
     \end{pmatrix},
     \begin{pmatrix}
    0&0 &0&a & b & c\\
     &a&b&c&0&0\\
  &&0&0&0&0\\
     &&&0&0&0\\
   &&&&0&0\\
    &&&&&0\\
     \end{pmatrix}.
\end{equation}
It appears that all these examples are constructed using 
 the \lq\lq\thinspace building blocks'' coming from the matrices (\ref{primo}) and (\ref{secondo}).
Other similar examples of spaces of dimension $\geq 4$  can be constructed using (\ref{secondo}).


\begin{rem} \label{rank1}
{\bf Kernel of rank $1$.} 

Note that for all $r\geq 3$, there are examples of $3$-dimensional vector spaces of matrices of order $2r+1$ and rank $2r$, corresponding to $\cK=\cO_{\p^{2}}(-r)$, generalizing (\ref{O2}). 

But $\dim SL(2r+1)<\dim\GG(2,\PP(\Lambda^2\kk^{2r+1}))$. Hence there are infinitely many orbits corresponding to the same kernel bundle.\end{rem}

\section{Vector spaces of skew--symmetric matrices of dimension three} \label{piani}

In this section we consider $3$-dimensional vector spaces $M$  of skew--symmetric matrices of order $8$ and constant rank $6$. 
Since $M$ has constant rank  $6$ then the vector bundle  $\cK$ in (\ref{fibrativettassociati}) has rank  $2$,  $c_1(\cK^{*})=3$ and  $\cK^{*}$  gives a $3$-Veronese embedding of $\PP(M)$ in $\GG(1,7)$, see Remark \ref{-K da gaussmap}.   Triple Veronese embeddings of $\PP^{n}$ in Grassmannians can be classified. 
\begin{thm}
\label{huh} Let $X\subset \GG(1,N)$ be a triple Veronese embedding of $\PP^{n}$  given by a vector bundle $E$ of rank $2$ on $\PP^{n}$  together with an epimorphism $\cO_{\PP^{n}}^{N+1}\to E$. Then one of the following holds:
\begin{itemize}
\item [(1)] $E \cong  \cO_{\PP^{n}}(a)\oplus \cO_{\PP^{n}}(3-a)$, $a=0,1$; $c_2(E)=0$ if $a=0$ and  $c_2(E)= 2$ if $a=1$;
\item [(2)]  $n=2$ and $E\cong \Omega_{\PP^{2}}(3)\cong T_{\PP^{2}}$,  $c_2(E)=3$; 
\item [(3)] $n=2$ and $E$ admits a resolution,
$$0\to\cO_{\PP^{2}}(2)\to E\to \cI_{p}(1)\to 0$$ 
where $\cI_{p}$ is the ideal sheaf of a point $p\in \PP^{2}$;  $c_2(E)=3$;
\item [(4)]  $n=2$ and $E$ is a stable vector bundle of rank $2$ on $\PP^{2}$ admitting one of the following resolutions 
\begin{itemize}
\item[(a)]
$0\to\cO_{\PP^{2}}(-1)^{\oplus 3}\to \cO_{\PP^{2}}^{\oplus 5} \to E\to  0$
\item[(b)]
$0\to \Omega_{\PP^{2}}(1)\oplus \cO_{\PP^{2}}(-1)^{\oplus 2} \to \cO_{\PP^{2}}^{\oplus 6} \to E\to  0$
\item[(c)]
$0\to \Omega_{\PP^{2}}(1)^{\oplus 2} \oplus \cO_{\PP^{2}}(-1)\to \cO_{\PP^{2}}^{\oplus 7} \to E\to  0$
\end{itemize}
\end{itemize}
In these last three cases $c_2(E)=6$ in $(a)$,  $c_2(E)=5$  in $(b)$,  $c_2(E)=4$   in $(c)$. 
\end{thm}

This result is due to S. Huh,  \cite[Theorem 1.1]{huh}, but in his theorem appears also another globally generated vector bundle over $\PP^{2}$ which  does not give an embedding of $\PP^{2}$ in $ \GG(1,N)$,  as J. C. Sierra has pointed out to us. More precisely the following lemma holds.

   \begin{lem} Let $E$ be a stable vector bundle on $\PP^{2}$ admitting the  resolution 
$ 0\to \Omega_{\PP^{2}}(1)\oplus \cO_{\PP^{2}}(-2)\to \cO_{\PP^{2}}^{\oplus 5} \to E\to  0$,
which corrisponds to the case  $(4)$, $(b)$ in \cite{huh}. $E$ does not give an embedding of $\PP^2$ in  $\GG(1, 4)$. 
 \end{lem}
 \proof    We can write the resolution of $E$ in the form
  \begin{eqnarray}
 \label{resol 4b}
0\to \Omega_{\PP^{2}}(1)\oplus \cO_{\PP^{2}}(-2)\to V\otimes \cO_{\PP^{2}} \to E\to  0
 \end{eqnarray}  
where $\dim V=5$.  The epimorphism $ V\otimes \cO_{\PP^{2}} \to E\to  0$ determines  a regular morphism  $\varphi_{V}: \PP^{2} \to \GG(1,4)$. 
Dualizing (\ref{resol 4b})  we get
 $$0\to E^{*} \to V^{*}\otimes \cO_{\PP^{2}} \to T_{\PP^{2}}(-1)\oplus \cO_{\PP^{2}}(2)\to   0.$$
The pair $(T_{\PP^{2}}(-1)\oplus \cO_{\PP^{2}}(2), V^{*})$ gives a triple Veronese embedding $\phi_{V^{*}}: \PP^{2} \to \GG(2,4)$ and $\varphi_{V}(\PP^{2})\cong \phi_{V^{*}}(\PP^{2})$ (see \cite[\S 4]{su}).  
Let $\overline{Y}$ be the $3$-fold in $\PP^{4}$ union of the lines of $\varphi_{V}(\PP^{2})$ (see  Facts \ref{facts}). Note also that the vector bundle  $T_{\PP^{2}}(-1)\oplus \cO_{\PP^{2}}(2)$ and  $V'=H^{0}(T_{\PP^{2}}(-1)\oplus \cO_{\PP^{2}}(2))$ give a triple Veronese embedding   $\varphi'_{V'}: {\PP^2} \to \GG(2,8)$, since  $h^{0}(T_{\PP^{2}}(-1)\oplus \cO_{\PP^{2}}(2))=9$. Let 
\begin{eqnarray*}
0 \to \cK' \to V'\otimes \cO_{\PP^2}\to T_{\PP^2}(-1)\oplus\cO_{\PP^2}(2)\to 0
\end{eqnarray*}
be the exact sequence of vector bundles associated to 
the epimorphism $V'\otimes \cO_{\PP^2}\to T_{\PP^2}(-1)\oplus\cO_{\PP^2}(2)\to 0$, so $V'^{*}\otimes \cO_{\PP^2}\to \cK'^{*}\to 0$ is also an epimorphism and hence the pair  $(\cK'^{*}, V'^{*})$ defines a map $\phi'_{V'^{*}}: {\PP^2} \to \GG(5,8)$ and $\phi'_{V'^{*}}({\PP^2})\cong \varphi'_{V'}({\PP^2})$ in the duality between $\GG(5,8)$ and $\GG(2,8)$. 
Let $\overline{Y}'$, $\overline{Z}'$ be the subvarieties of ${\PP^8}$ associated to $\phi'_{V'^{*}}({\PP^2})$ and $\varphi'_{V'}({\PP^2})$, respectively. The $3$-fold  $\overline{Y}\subset {\PP^4}$ is obtained after slicing $\overline{Y}'\subset {\PP^8}$ with $4$ hyperplanes; this passes from $\phi'_{V'^{*}}({\PP^2})\subset  \GG(5,8)$ to  $\varphi_{V}({\PP^2}) \subset  \GG(1,4)$. By duality this is equivalent to projecting $\varphi'_{V'}({\PP^2})\subset \GG(2,8)$ in $\GG(2,4)$ and successively dualizing to $\GG(1,4)$. Since $\varphi'_{V'}$ is given by $(T_{\PP^2}(-1)\oplus \cO_{\PP^2}(2), V')$, the variety $\overline{Z}'$ corresponding to  $\varphi'_{V'}({\PP^2})$ is constructed as follows: fix in ${\PP^8}$ a ${\PP^2}$ and a ${\PP^5}$ complementary to each other and a $v_2({\PP^2})$ in ${\PP^5}$, fix an isomorphism $\psi$ between ${\PP^2}^{*}$ and $v_2({\PP^2})$ and consider the family of the $2$-planes joining a line and a point corresponding to each other in $\psi$. When we project in ${\PP^4}$, we get a 
 ${\PP^2}=\alpha$ and a projected Veronese surface, of degree $4$, intersecting $\alpha$ in $4$ points. Hence $4$ planes of the family come together to coincide with $\alpha$. This gives rise to a  point of multiplicity $4$ of $\phi_{V^{*}}({\PP^2})$. Hence also $\varphi_{V}({\PP^2})$ is singular.
 \qed
 \vspace{2mm}

\begin{rem} To explain how the classification in Theorem \ref{huh} is organized, we note that the bundles in (1), (2) are uniform, while those in (3), (4) are not. Moreover the bundle appearing in (3) is unstable, while those in (2) and (4) are stable. The corresponding moduli spaces $\mathcal M(3, c_2)$ have dimension $4c_2-12$ (see \cite[Ch. 2, \S 4]{OSS}). 
\end{rem}


We will see that the non split vector bundles given in Theorem \ref{huh}  can be seen as  quotient of vector bundles of higher rank of a very particular form. This fact turns out to be  crucial in constructing skew--symmetric matrices of constant rank $6$. 

\begin{definition}(\cite{su}) We say that a vector bundle $\cF$ on ${\PP^n}$ is a quotient of $\cE$ if there exists an exact sequence 
$0\to \cO_{\PP^n}^{\oplus s} \to \cE \to \cF \to 0$, 
corresponding to $s$ sections of $\cE$. 
\end{definition}

\begin{prop} 
\label{quotientof}
Let $E$ be a vector bundle of rank $2$ on ${\PP^2}$ defining a triple Veronese embedding of ${\PP^2}$ in a Grassmannian $\GG(1,N)$, as in Theorem \ref{huh}. 
\begin{itemize}
\item[(i)] If $E$ is as in $(2)$ then $E$ is a quotient of  $\cO_{\PP^2}(1)^{\oplus 3}$;
\item[(ii)]  If $E$ is as in $(3)$ then $E$ is a quotient of  $T_{\PP^2}(-1)\oplus \cO_{\PP^2}(2)$;
\item[(iii)]  If $E$ is  as in $(4), (a)$ then $E$ is a quotient of $T_{\PP^2}(-1)^{\oplus 3}$;
\item[(iv)] If $E$ is  as in $(4), (b)$ then $E$ is a quotient of  $T_{\PP^2}(-1)^{\oplus 2}\oplus \cO_{\PP^2}(1)$;
\item[(v)] If $E$ is  as in $(4), (c)$ then $E$ is a quotient of $T_{\PP^2}(-1)\oplus \cO_{\PP^2}(1)^{\oplus 2}$.
\end{itemize}
\end{prop}
\proof   The case $(i)$ follows from the Euler exact sequence.

In the case $(ii)$ the vector bundle $E$  is unstable because $h^{0}(E_{norm})=h^{0}(E(-2))\neq 0$. Let  $\cG=T_{\PP^2}(-1)\oplus \cO_{\PP^2}(2)$, then   $c_1(\cG)=c_2(\cG)=3$,   moreover  $h^{0}(\cG)\neq 0$ hence  there exists an exact sequence
\begin{eqnarray}
\label{quoziente}
0\to \cO_{\PP^2}\to \cG \to \cQ\to 0
\end{eqnarray}
 corresponding to a section of $\cG$. 
Note that $\cQ$ normalized, $\cQ_{norm}=\cQ(-2)$. Twisting (\ref{quoziente}) with $\cO_{\PP^2}(-2)$ and considering its associated  cohomology exact sequence  it follows that $h^{0}(\cQ(-2))=1$ and thus $\cQ$ cannot be stable, see \cite[Lemma 1.2.5, pg 165]{OSS}. Moreover, because $c_1(\cQ)=3$ and $rk(\cQ)=2$, then by 
  \cite[ Remark 1.2.3 pg 163]{OSS} it follows that $\cQ$ is stable if and only if is semistable. Hence, being $h^{0}(\cQ(-2))=1$,  we can conclude that  $\cQ$ is unstable and thus $\cQ$ has to be the vector bundle $E$ in $(2)$, since $c_2(\cQ)=3$.\\
  In the case $(iii)$  the resolution of $E$ along with the Euler exact sequence yields the following commutative diagram with exact rows and columns
\begin{eqnarray*} \label{4-a}
    \xymatrix{ & &0\ar^-{}[d] &0\ar^-{}[d] &\\
    & &\cO_{\PP^2}^{\oplus 4}  \ar^-{}[d] \ar@{=}[r]  &\cO_{\PP^2}^{\oplus 4}  \ar^-{}[d]&\\
    0 \ar[r] & \cO_{\PP^2}(-1)^{\oplus 3}\ar@{=}[d] \ar[r] &  \cO_{\PP^2}^{\oplus 9}
      \ar^-{}[d] \ar[r] & 
      T_{\PP^2}(-1)^{\oplus 3}\ar[r] \ar^-{}[d] &  0 \\
      0 \ar[r] & \cO_{\PP^2}(-1)^{\oplus 3} \ar[r] & \cO_{\PP^2}^{\oplus 5}\ar^-{}[d]  \ar[r] & E \ar^-{}[d] \ar[r]  & 0 \\
    & &0&0&}
    \end{eqnarray*}
and we get that $E$ is a quotient of $T_{\PP^2}(-1)^{\oplus 3}$.\\
 The proof of $(iv)$ and $(v)$ runs along the same lines of  $(iii)$, hence we omit it.   \qed
 \vspace{2mm}
 
\begin{corol}\label{corembed}  Let $E$ be a rank $2$ vector bundle on $\PP^2$ with $c_1(E)=3$. If $E$ gives an embedding of $\PP^{2}$ in $\GG(1,N)$  then either $E$ splits or  $E$ is a quotient of a vector bundle of higher rank which is a direct sum of $T_{\PP^2}(-1)$(or more copies of it)  and  $\cO_{\PP^{2}}(k)$ (or more copies of it) for some positive integer $k$.
\end{corol}

\begin{rem} Let $M$ be a vector space of skew--symmetric matrices of constant rank $2r$ and order $N+1$. We can interpret it as a linear space contained in $S_{r-1}\GG(1,N)\setminus S_{r-2}\GG(1,N).$ Observe that the exact sequence (\ref{fibrativettassociati1}), after taking a quotient $\cQ$ of $\cK$, gives rise to another exact sequence in which a new matrix comes up and one wants to  know if its rank is constant. Taking a quotient $\cQ$ of $\cK$  corresponds to projecting $\PP^N$ to $\PP^{N-1}$ from a point $O$.
The projection $\pi_O$ of centre $O$   induces a projection $\pi_{\Lambda_O}$ from $\PP(\Lambda^2\kk^{N+1})$ to $\PP(\Lambda^2\kk^{N})$, whose centre is the subspace $\Lambda_O\subset \GG(1,N)$, representing the lines through $O$. How should one choose the centre of projection  in order that the rank of $M$ remains constant under this projection? The answer is given by the following Proposition.
\end{rem}

We recall that a point $\omega$ in $S_{r-1}\GG(1,N)\setminus S_{r-2}\GG(1,N)$ can be written in the form $[v_1\wedge w_1+\ldots +v_r\wedge w_r]$, where $v_1,\ldots,v_r,w_1,\ldots,w_r$ are linearly independent vectors; the corresponding points generate a subspace $L_\omega$ of $\PP^n$ of dimension $2r-1$. Then the entry locus of $\omega$ is the subgrassmannian  $\GG(1,L_\omega)$,  namely a point of $\GG(1,N)$ belongs to some $(r-1)$-plane $r$-secant to $\GG(1,N)$ and containing $\omega$ if and only if it belongs to $\GG(1,L_\omega)$.

\begin{prop}
\label{projection}
Let $\PP(M)\subset S_{r-1}\GG(1,N)$ be a vector space of matrices of constant rank $2r$. Let $O\in \PP^N$ be a point such that $\PP(M)\cap \Lambda_O=\emptyset$. Then the matrices of $\pi_{\Lambda_O}(\PP(M))$ have  constant rank $2r$ if and only if $O$ does not belong to the union of the spaces $L_\omega$, as $\omega$ varies in  $\PP(M)$.
\end{prop}
\proof  Let $\omega=[v_1\wedge w_1+\ldots +v_r\wedge w_r]$ be a point of $\PP(M)$. Then  $\pi_{\Lambda_O}(\PP(M))(\omega)=[Av_1\wedge Aw_1+\ldots +Av_r\wedge Aw_r]$, where $A$ is a matrix representing $\pi_O$, and its rank  is strictly less than $r$ if and only if $v_1,\ldots,v_r,w_1,\ldots,w_r$ can be chosen so that some summand $Av_i\wedge Aw_i$ vanishes. But this means precisely that $O$ belongs to $L_\omega$.
\qed
\vspace{2mm}

\begin{corol}\label{projmat} Let $\PP(M)$ be a linear space of matrices of constant rank $2r$ and dimension $d$. Then $\PP(M)$ can be isomorphically projected to $S_{r-1}\GG(1,2r+d-1)$ so that its rank remains constant and equal to $2r$.
\end{corol}
\proof It is enough to note that $\dim \bigcup_{\omega\in \PP(M)}L_\omega\leq \dim \PP(M)+2r-1.$\qed
\vspace{2mm}

In particular a projective $2$-plane of matrices of constant rank $6$ can be projected in $S_2\GG(1,7)$ mantaining constant rank $6$.

We can now state the main result of this section, which gives a reverse statement to Proposition \ref{-K da gaussmap}.
\begin{thm}\label{proj}
Let $E$ be a rank two vector bundle on $\PP^2$ defining a triple Veronese embedding of $\PP^2$ in $\GG(1,7)$. Then there exists a vector space of $8\times 8$ matrices of constant rank $6$ whose associated bundle $\cK$ is such that $E\simeq\cK^*$.
\end{thm} 
\proof By Corollary \ref{corembed} $E$ is  a direct sum of copies of $\cO_{\PP^2}$, $\cO_{\PP^2}(1)$, $\cO_{\PP^2}(2)$, $\cO_{\PP^2}(3)$, $T_{\PP^2}(-1)$, or a quotient of it. From the results of Section \ref{dim maggiore}, each of these bundles is the dual of a bundle appearing as kernel in an exact sequence of the form (\ref{fibrativettassociati}).  Taking a direct sum of   matrices corresponding to the direct summands, we construct a matrix of constant rank $6$ and  order possibly bigger than $8$. Finally, by Corollary \ref{projmat}, with a suitable projection we get a $8\times 8$ matrix of the desired form.
\qed
\vspace{2mm}

\subsection{Examples} For each class of rank two bundles appearing in Theorem \ref{huh}, we will give now one or more examples of linear systems of matrices of constant rank $6$. Unfortunately we are not able to give a complete classification of the orbits for the action of $SL(8)$. As we have already noted in Subsection \ref{rank1}, for some bundles there are infinitely many orbits. On the other hand, for the bundles in Theorem \ref{huh}, (4), there is a moduli space of positive dimension.
\begin{exm}\label{split}  {\bf Split bundles.}

 Let $\pi_{1}$ be the plane 
\begin{eqnarray*}
\begin{pmatrix}
        0 & 0 & 0 & 0 & a & b & 0\\
          & 0 & 0 & a & b & c & 0\\
       &  & a & b & c & 0 & 0 \\
       & &  & c & 0 & 0 & 0\\
         & & & & 0 & 0 & 0\\
        &  && &  & 0 & 0\\
         & &  &  &  & & 0\\
       \end{pmatrix}.
\end{eqnarray*}
            
\noindent  The associated rank $2$ vector bundle  is  $\cK^{*}=\cO_{\PP^2}\oplus\cO_{\PP^2}(3)$, $c_2(\cK^{*})=0.$

Let $\pi_{2}$ be the plane 
\begin{eqnarray*} \begin{pmatrix}
          0 & 0 & a & b & 0 & 0 & 0\\
        & a & b & c & 0 & 0 & 0\\
        && c & 0 & 0 & 0 & 0 \\
        &  & & 0 & 0 & 0 & 0\\
         &  &&  & 0 & 0 & 0\\
       &&&& & a & b\\
          &  & && &  & c\\
       \end{pmatrix}.
        \end{eqnarray*}
          \end{exm}
    \noindent     The associated rank $2$ vector bundle  is  $\cK^{*}=\cO_{\PP^2}(1)\oplus\cO_{\PP^2}(2)$, $c_2(\cK^{*})=2.$
          
           \begin{exm} \label{Steiner}  {\bf Steiner bundles}
           
     \noindent Steiner bundles are quotients of $T_{\PP^2}(-1)^3$, have $c_2=6$ and move in a moduli space of dimension $12$. An example of matrix is   
     \begin{eqnarray*} \begin{pmatrix}\label{schwarz}
       0 & 0 & a & b & c & 0 & 0\\
         & 0 & 0 & a & b& c & 0\\
    && 0 & 0 & a & b & c \\
        &&  & 0 & 0 & 0 & 0\\
        &&  &  & 0 & 0 & 0\\
  & &  & &  & 0 & 0\\
         & &  & &  & & 0\\
       \end{pmatrix}
       \end{eqnarray*}        
\noindent       In this case the associated  bundle   $\cK^{*}$ is a Schwarzenberger bundle,  having a conic of jumping lines. General Steiner bundles have $6$ jumping lines. Following the construction of Dolgachev-Kapranov (see \cite{dkDuke}) and choosing  as follows the equations of the jumping lines:

$x_0=0;$

$x_1=0;$

$x_2=0;$

$\lambda_0x_0+\lambda_1x_1+\lambda_2x_2=0;$

$\mu_0x_0+\mu_1x_1+\mu_2x_2=0;$

$\nu_0x_0+\nu_1x_1+\nu_2x_2=0,$

\noindent we get a matrix $M$ of the form:
\begin{equation} M=
\begin{pmatrix}
         0 & B \\
        - B^t & 0 \\
 \end{pmatrix}
 \end{equation} 
 where 
\begin{equation}B=
\begin{pmatrix}
         \lambda_0a+\mu_0b+\nu_0c & 0 & 0 & \lambda_0a & \mu_0b \\
         0 & \lambda_1a+\mu_1b+\nu_1c & 0 & \lambda_1a & \mu_1b \\
         0 & 0 & \lambda_2a+\mu_2b+\nu_2c & \lambda_2a & \mu_2b \\
         \end{pmatrix}.
 \end{equation}

Hence we have  a family of examples, depending on the parameters $\lambda, \mu,\nu$.  Such parameters have to be chosen so that the six jumping lines are in general position. 
      \end{exm}

        \begin{exm}\label{unstable}   {\bf Unstable bundle}
       
       Let $\pi_{3}$ be the plane 
       \begin{eqnarray*}\begin{pmatrix}
        0 & 0 & 0 & 0 & a & b & c\\
  & 0 & 0 & a & b& 0 & 0\\
         &  & a & b & c & 0 & 0 \\
    &&  & c & 0 & 0 & 0\\
        & & &  & 0 & 0 & 0\\
         & &  & & & 0 & 0\\
          & &  & &  & & 0\\
       \end{pmatrix}
   \end{eqnarray*}
       \end{exm}

\noindent  In this case the associated rank $2$ vector bundle   $\cK^{*}$ is the unstable one,  quotient of  $T_{\PP^2}(-1)\oplus\cO_{\PP^2}(2)$, $c_2(\cK^{*})=3$ (case (3) of Theorem \ref{huh}).

  \begin{exm}\label{ex4}    
       Let $\pi_{4}$ be the plane 
    \begin{eqnarray*}\begin{pmatrix}
        0 & 0 & 0 & 0 & a & b & c\\
          & 0 & 0 & a & b& c & 0\\
         && a & b & 0 & 0 & 0 \\
          & &  & c & 0 & 0 & 0\\
         &  &&  & 0 & 0 & 0\\
      & &  & &  & 0 & 0\\
& &  & &  &    & 0\\
       \end{pmatrix} 
\end{eqnarray*}
 \end{exm}
\noindent  In this case the associated rank $2$ vector bundle   $\cK^{*}$ is a quotient of  $T_{\PP^2}(-1)\oplus T_{\PP^2}(-1)\oplus\cO_{\PP^2}(1)$, $c_2(\cK^{*})=5$ (case (4)(b) of Theorem \ref{huh}).

\vspace{2mm}
The expression of the matrices in the following Examples \ref{ex5} and  \ref{ex6}  is not so evident a priori.  The matrices in these examples correspond, respectively,  to the bundle (4)(c) of Theorem \ref{huh}  and to the tangent bundle $T_{\PP^2}$. 

  \begin{exm}\label{ex5}    {\bf A quotient of $T_{\PP^2}(-1)\oplus\cO_{\PP^2}(1)\oplus\cO_{\PP^2}(1)$.}
       
       Projecting from $[1,0,0,0,1,0,0,0,0,1] $ and subsequently from $[0,0,1,1,0,0,0,1,0]$ the direct sum matrix naturally associated to this bundle, we get the following plane   $\pi_{5}$: 
       \begin{eqnarray*} \begin{pmatrix}
         c & a & 0 & 0 & 0 & 0 & a\\
         & b & 0 & 0 & 0 & 0& b\\
          & & c-b & 0 & 0 & a & 0 \\
         & & & 0 & 0 & b & 0\\
          & &  &  & a & b & c\\
 & & &&  & 0 &0\\
         & &&  &&  & 0\\
                \end{pmatrix}
     \end{eqnarray*}
       \end{exm}

 \begin{exm}\label{ex6}    {\bf The tangent bundle.}
      
      Let $\pi_{6}$ be the plane 
      \begin{eqnarray*}\begin{pmatrix}
         c & 0 & a & 0 & 0 & 0 & a\\
        & 0 & b & 0 & 0 & 0& b\\
           &  & a & b & 0 & 0 & 0 \\
      & & & c & 0 & 0 & 0\\
         & & & & 0 & 0 & 0\\
         && & &  & a &b\\
          &  &  &&&  & c\\
       \end{pmatrix}
 \end{eqnarray*}
 \end{exm}
\noindent  It has been constructed projecting from $[1,0,0,0,1,0,0,0,1]$ the direct sum matrix, coming from the expression of    $\cK^{*}=T_{\PP^2}$, as a quotient of $\cO_{\PP^2}(1)\oplus\cO_{\PP^2}(1)\oplus\cO_{\PP^2}(1)$.

       \subsection{Triple Veronese embeddings.}

Let $\gamma$ be the Gauss map from $S_2\GG(1,7)$ to $\GG(1,7)$. We will shortly give now the geometrical interpretation of the varieties $\gamma(\PP(M))$ for each of the above examples. They are all (projections of) the Veronese variety $v_3(\PP^2)$.

\begin{itemize}
   \item  $\gamma(\pi_1)$ is contained in $\GG(1,6)$ and represents the lines of a cone with vertex one point over 
  $v_3(\PP^2)$ projected from $\PP^9$ to $\PP^6$. Since varying the centre of projection we get varieties isomorphic but not always projectively equivalent, this explains the presence of infinitely many orbits for these planes. 

\item $\gamma(\pi_2)\subset \GG(1,7)$  represents the lines joining the corresponding points in an isomorphism between  a fixed $2$-plane and a projected $2$-Veronese surface in a fixed $\PP^4$.

\item If $M$ is one of the planes of Example \ref{Steiner}, then $\gamma(\PP(M))$ is contained in a subgrassmannian $\GG(1,4)$, where the $\PP^4$ is defined by the equations $x_0=x_1=x_2=0$. It represents the lines of a cubic threefold, whose equation is the determinant of a $1$-generic matrix $\Omega$ obtained as follows. Write $B=aB_1+bB_2+cB_3$, and let $B_i^j$ denote the $j$-th row of $B_i$. Let $X'$ be the transposed of $(x_3 \ \ldots \ x_7)$. Then 
\begin{eqnarray}\Omega=\begin{pmatrix}
B_1^1X'&B_1^2X'&B_1^3X'\\
B_2^1X'&B_2^2X'&B_2^3X'\\
B_3^1X'&B_3^2X'&B_3^3X'\\
\end{pmatrix}\end{eqnarray}
In the special case of the Schwarzenberger bundle, we get the cubic threefold of the secant lines of a quartic rational normal curve.

\item The lines of $\gamma(\pi_3)$ are obtained as follows. Note that, since 
 $\dim H^{0}(\PP^2, T_{\PP^2}(-1)\oplus \cO_{\PP^2}(2))=9$, this bundle  gives a triple Veronese embedding of $\PP^2$ in $\GG(2, 8)$. Geometrically we  fix an isomorphism between $v_2(\PP^2)$ and ${\overset\vee{\ \PP^2}}$  and we get a family of $\PP^{2}$'s spanned by a point in $v_2(\PP^2)$ and the corresponding line in ${\overset\vee{\ \PP^2}}$.  Taking a quotient of  $T_{\PP^2}(-1)\oplus \cO_{\PP^2}(2)$ is the same as cutting this family with a hyperplane,  and this  gives $\gamma(\pi_3)\subset \GG(1, 7)$.  

 \item The description of $\gamma(\pi_4)$ and $\gamma(\pi_5)$ is similar to the previous one.
Since $\dim H^{0}(\PP^2, T_{\PP^2}(-1)\oplus T_{\PP^2}(-1) \oplus \cO_{\PP^2}(1))=9$ we have a triple Veronese embedding of $\PP^2$ in $\GG(4, 8)$. Geometrically we fix three planes and we have a corrispondence between  the first plane $\PP^2$ and the dual of the other two ${\PP^2}$'s.  We get a family of $\PP^{4}$'s spanned by a point in $\PP^2$ and the two corresponding lines in the two ${\overset\vee{\ \PP^2}}$.  Cutting this family with three hyperplane we get  $\gamma(\pi_4)$.  If we consider instead the bundle $\cO_{\PP^2}(1)^{\oplus 3}$ this gives an triple Veronese embedding of $\PP^2$ in $\GG(2, 8)$.  We consider again three disjoint  $\PP^2$'s and we get the family of  $\PP^2$'s spanned by three corresponding points in fixed isomorphisms among them. Cutting this family with a hyperplane we get  $\gamma(\pi_5)$.   
   \end{itemize}

\begin{rem} The algorithm  in Section \ref{algo} can be suitably modified to compute the dimensions of the orbits of the matrices constructed in this Section. One obtains
that the dimension of the orbits  is $54$ (respectively $60$) in    
Example \ref{split}; $52$ (respectively $56$)  in   
Example \ref{Steiner};  $58$ in  Example \ref{unstable}  and in Example \ref{ex4}; $59$ in Example \ref{ex5} and $60$ in Example \ref{ex6}.  \end{rem}

\section{Westwick example rivisited} \label{dim3}
We start this final section with the following:
\begin{rem}
There do not exist vector spaces of dimension $4$ of $8\times 8$ matrices of constant rank $6$. This follows from a computation on the Chern classes, see for instance  \cite[Example 2.12]{IL}.
\end{rem}

The first possibility  for  a $\PP^3$ of skew--symmetric matrices of order $2r+2$ and constant rank $2r$ is for $r=4$. The only known example has been given by Westwick (\cite{w2}) and is  the following: 

 \begin{equation}\label{west}
     \begin{pmatrix}
      0&0 &0& 0&0& 0&a & b & 0\\
     &0& 0&0& 0&a&b&0&c\\
   && 0&0&-a&b&0&c&d\\
     &&&a&b&0&c&d&0\\
     &&&&0&c&-d&0&0\\
     &&&&&d&0&0&0\\
     &&&&&&0&0&0\\
    &&&&&&&0&0\\
     &&&&&&&&0\\
     \end{pmatrix}
\end{equation}

We will say something about the vector bundles associated to such $\PP^3=\PP(M)$ of skew--symmetric matrices of constant rank $8$ and order $10$.  

With the notation in Section 1, one can easily see, using for instance \cite[Example 2.12]{IL}, that $c_{1}(\cK^*)=4$ and $c_{2}(\cK^*)=6$. Let $s\in H^{0}(\PP^{3}, \cK^*)$ be a generic section and let $Y$ be  its scheme of zeros. Because $\cK^*$ is spanned by global sections then $Y$ is smooth of codimension $2=rk (\cK^*)$. The section $s$ defines an exact sequence 
 \begin{equation}
 \label{successione associata a s}
 0\to  \cO_{\PP^{3}} \to \cK^{*} \to {\mathcal J}_{Y}(c_{1}(\cK^*))\to 0
 \end{equation}
We see that $deg(Y)=c_2(\cK^*)=6$, $N_{Y/{\PP^{3}}}= \cK^*_{|Y}$.  We get, by adjunction,  that  $K_{Y}= \cO_{Y}$ and thus $g(Y)=1.$ 
 Twisting the exact sequence (\ref{successione associata a s}) with $\cO_{\PP^{3}}(-2)$ and recalling that $c_{1}(\cK^*)=4$ we get
  \begin{equation}
 \label{successionetwistata}
 0\to  \cO_{\PP^{3}}(-2) \to \cK^*(-2) \to {\mathcal J}_{Y}(2) \to 0
 \end{equation}
From the cohomology sequence associated to (\ref{successionetwistata}), using the fact that $H^{0}(\PP^{3},  {\mathcal J}_{Y}(2))=0$ because $Y$ cannot be contained in any  quadric surface, it follows that $H^{0}(\PP^{3}, \cK^*(-2))=H^{0}(\PP^{3}, \cK^*_{norm})=0$ and thus $\cK^*$ is a stable vector bundle. 

By Proposition \ref{gaussmap} we see that $\gamma(\PP(M))$ is  a $4$-tuple Veronese embedding of $\PP^{3}$ in $\GG(1,9)$. Note that this embedding is given by a proper subspace of $H^{0}(\PP^{3},   \cK^*)$. In fact, using results contained in \cite{hartshorne}, one can show  that $h^{0}(\PP^{3},   \cK^*)=12$.

Thus on $\PP^{3}$ such vector bundles $\cK^{*}$ are the only  ones for which it can exists a $10\times 10$ skew--symmetric matrix of constant rank $8$.  R. Hartshorne in \cite[Corollary 9.8]{hartshorne} has proved that  the variety of moduli of these bundles is an irreducible nonsingular variety of dimension 
$13$.
\begin{rem} From (\ref {successione associata a s}) one   computes also that $h^{1}(\PP^{3},   \cK^*)=0$. Hence $\cK^*$ is not a quotient of any bundle of higher rank. Similarly, the restriction $\cK^*\mid_H$ of $\cK^*$ to a general plane $H$ is a stable bundle with $h^0(H,\cK^*\mid_H)=10$ and $h^1(\cK^*\mid_H)=0$. Therefore the matrix obtained by restricting  (\ref{west}) to a general plane can be thought of as a new building block for constructing vector spaces of dimension $3$ of matrices of constant rank $\geq 8$ and order at least $10$.
\end{rem}

\begin{rem} The study of linear spaces of skew--symmetric matrices of constant rank is  related to the study of possible degenerations of a class of projective varieties called {\it Palatini scrolls}, that is those varieties $X$ in $\PP^{N}$, with $N$ odd, which are
degeneracy loci of general morphisms $\phi:
\cO^{m}_{\PP^{N}} \to \Omega_{\PP^{N}}(2)$. Such $X$ is smooth if $m<\frac{N+4}{2}$, \cite{bazan-mezzetti}. As it is well known a morphism $\phi: \cO^{m}_{\PP^{N}} \to \Omega_{\PP^{N}}(2)$ gives a $(N+1)\times (N+1)$ skew--symmetric matrix  of linear forms $M_{\phi}$ on ${\PP^{m-1}}$.

For instance, if $N=4$ and $m=3$, then $X$ is a projected Veronese surface, if $N=5$ and $m=4$, then $X$ is a Palatini threefold: its degenerations have been studied in \cite{dm} relying on the classification given in \cite{MM}.

If $m=5$ and $N=7$ then $X$ is a smooth fourfold in ${\PP^{7}}$ with base of the scroll the quartic $3$-fold $Y$ in ${\PP^{m-1}}=\PP^4$, defined by $\Pf (M_{\phi})$. The fact that there do not exist vector spaces of dimension $4$ of $8\times 8$ matrices of constant rank $6$ says  that there cannot be a degeneration of $X$ obtained by degenerating the base $Y$ so that  it acquires a $\PP^3$ as irreducible component.

The situation is different for $m=5$ and $N=9$. In this case the base $Y$ of the scroll is a quintic $3$-fold  defined by $\Pf (M_{\phi})$, with $M_{\phi}$ a $10\times 10$ matrix of linear forms on $\PP^4$.
From the example  (\ref{west})  it follows that  the base $Y$ can degenerate so that  it contains a $\PP^3$ as irreducible component.
\end{rem}

\section*{Acknowledgements} 
We would like to thank Jos\'e Carlos Sierra  for many interesting discussions and suggestions.

We also like to thank Giorgio Ottaviani for  suggesting us the Macaulay  script   in \ref{algo} which helps in computing  the dimension of the orbits.


\end{document}